\documentclass[twosided,12pt,reqno,]{article}
\usepackage{amssymb}
\usepackage{amsthm}
\usepackage{amsmath}
\usepackage{a4wide}
\usepackage[all]{xy}
\usepackage{fancyhdr}

\newtheorem{theorem}{Theorem}[section]
\newtheorem{proposition}[theorem]{Proposition}
\newtheorem{corollary}[theorem]{Corollary}
\newtheorem{lemma}[theorem]{Lemma}
\theoremstyle{definition}
\newtheorem*{notation}{Notation}
\newtheorem*{Beweis}{Proof}
\newtheorem{definition}[theorem]{Definition}
\newtheorem{punto}[theorem]{}
\theoremstyle{remark}
\newtheorem{remark}[theorem]{Remark}
\newtheorem{ex}[theorem]{Example}
\newtheorem{exs}[theorem]{Examples}
\newtheorem{c-ex}[theorem]{Counterexample}
\newtheorem{c-exs}[theorem]{Counterexamples}
\newtheorem{remarks}[theorem]{Remarks}
\swapnumbers
\CompileMatrices
\input{tcilatex}
\begin{document}

\title{A Dual Zariski Topology for Modules\thanks{%
MSC2010: 16N20; 16N80 (13C05, 13C13, 54B99) \newline
Keywords: Duo Module; Multiplication Module; Comultiplication Module; Fully
Coprime Module; Zariski Topology; Dual Zariski Topology}}
\author{\textbf{Jawad Abuhlail}\thanks{%
The author would like to acknowledge the support provided by the Deanship of
Scientific Research (DSR) at King Fahd University of Petroleum $\&$ Minerals
(KFUPM) for funding this work through project No. FT08009.} \\
Department of Mathematics and Statistics\\
King Fahd University of Petroleum $\&$ Minerals \\
abuhlail@kfupm.edu.sa}
\date{\today}
\maketitle

\begin{abstract}
We introduce a \textit{dual Zariski topology} on the spectrum of \textit{%
fully coprime} $R$-submodules of a given \textit{duo} module $M$ over an
associative (not necessarily commutative) ring $R$. This topology is defined
in a way dual to that of defining the Zariski topology on the prime spectrum
of $R$. We investigate this topology and clarify the interplay between the
properties of this space and the algebraic properties of the module under
consideration.
\end{abstract}

\section{Introduction}

\qquad Inspired by the interplay between the Zariski topology defined on the
prime spectrum of a commutative ring $R$ and the ring theoretic properties
of $R$ in \cite{Bou1998, AM1969} (see also \cite{LY2006, ST2010, ZTW2006}),
we introduce in this paper a dual Zariski topology on the spectrum of \emph{%
fully coprime }submodules of a given non-zero \emph{duo} module $M$ over an
associative ring $R$ and study the interplay between the properties of $_{R}M
$ and the topological space we obtain. The spectrum we consider was
introduced for bicomodules over corings in \cite{Abu2006} (see also \cite%
{Wij2006}, \cite{WW2009}) and topologized in \cite{Abu2008} where a
Zariski-like topology was investigated. Some of the results in this paper
can be considered as module-theoretic versions of results in \cite{Abu2006}
and are dual to results in \cite{Abu} which is devoted to a Zariski topology
on the spectrum of fully prime submodules of a given duo module.

This paper extends also the study of the so called \emph{top modules,} i.e.
modules whose spectrum of \emph{prime submodules} attains a Zariski topology
(e.g. \cite{Lu1984, Lu1999, MMS1997, MMS1998, Zha2006-a, Zha2006-b}), to a
notion of primeness not dealt with so far from the topological point of view
(other notions appear, for example, in \cite{Dau1978, Wis1996, RRRF-AS2002,
RRW2005, Wij2006, Abu2006, WW2009}).

The paper is organized as follows: After this introduction, we collect in
Section 2 some preliminaries and recall some properties and notions from
Module Theory that will be needed in the sequel. In Section 3, given a
non-zero duo left $R$-module $M,$ we introduce and investigate a dual
Zariski topology on the spectrum $\mathrm{Spec}^{\mathrm{fc}}(M)$ of
non-zero submodules that are \emph{fully coprime} in $M$ (e.g. \cite%
{Abu2006, Wij2006, WW2009}). In particular, we investigate when this space
is Noetherian or Artinian (Theorem \ref{noth-art}), irreducible (Theorem \ref%
{corad-fc}), ultraconnected (Proposition \ref{uniform}), compact or locally
compact (Theorem \ref{count-compact}), connected (Theorem \ref{colocal}), $%
T_{1}$ (Proposition \ref{T1}) or $T_{2}$ (Theorem \ref{T2}).

\section{Preliminaries}

\qquad In this section, we fix some notation and recall some definitions and
basic results. For any undefined terminology, the reader is referred to \cite%
{Wis1991} and \cite{AF1974}.

Throughout, $R$ is an associative (not necessarily commutative) ring with $%
1_{R}\neq 0_{R}$ and $M$ is a non-zero unital left $R$-module. By an ideal
we mean a \textit{two-sided ideal} and by an $R$-module we mean a \textit{%
left} $R$-module, unless explicitly otherwise mentioned. We set $S:=\mathrm{%
End}(_{R}M)^{op}$ (the ring of $R$-linear endomorphisms of $M$ with
multiplication given by the opposite composition of maps) and consider $M$
as an $(R,S)$-bimodule in the canonical way. We write $L\leq _{R}M$ ($%
L\lvertneqq _{R}M$) to indicate that $L$ is a (proper)\ $R$-submodule of $M$%
. For non-empty subsets $L\subseteq M$ and $I\subseteq R$ we set
\begin{equation*}
(L:_{R}M):=\{r\in R|\,rM\subseteq L\}\,\text{ and }\,(L:_{M}I):=\{m\in
M|\,Im\subseteq L\}.
\end{equation*}

\begin{definition}
We say $L\leq _{R}M$ is \emph{fully invariant} or \emph{characteristic} iff $%
f(L)\subseteq L$ for every $f\in S$ (equivalently iff $L\leq M$ is an $(R,S)$%
-subbimodule). In this case, we write $L\leq _{R}^{\mathrm{f.i.}}M.$ We call
$_{R}M$ \emph{duo} or \emph{invariant} iff every $R$-submodule of $M$ is
fully invariant.
\end{definition}

Recall that the ring $R$ is said to be \textit{left duo} (\textit{right duo}%
) iff every left (right) ideal of $R$ is two-sided and to be \textit{left
quasi-duo} (\textit{right quasi-duo}) iff every maximal left (right) ideal
of $R$ is two-sided. Moreover, $R$ is said to be (\textit{quasi-}) \textit{%
duo} iff $R$ is left and right (quasi-) duo.

\begin{punto}
\label{ex-duo}

Examples of duo modules are:

\begin{enumerate}
\item uniserial Artinian modules \cite{OHS2006}.

\item self-injective self-cogenerator modules with commutative endomorphism
rings (e.g. \cite[48.16]{Wis1991}).

\item \textit{multiplication modules}: $_{R}M$ is multiplication iff every $%
L\leq _{R}M$ is of the form $L=IM$ for some ideal $I$ of $R$, equivalently $%
L=(L:_{R}M)M$. Multiplication modules over commutative rings have been
studied intensively in the literature (e.g. \cite{AS2004, PC1995, Smi1994}).
Several results in these paper have been generalized by Tuganbaev (e.g. \cite%
{Tug2003}, \cite{Tug2004}) to modules over rings \textit{close to be
commutative} (see \cite{Abu} for a summary).

\item \textit{comultiplication modules}: $_{R}M$ is comultiplication iff
every $L\leq _{R}M$ is of the form $L=(0:_{M}I)$ for some ideal $I$ of $R$,
equivalently, $L=(0:_{M}(0:_{R}L))$. A commutative ring for which $_{R}R$ is
a comultiplication module is called a \emph{dual ring}. For examples and
results on such modules, we refer mainly to \cite{AS} and \cite{A-TF2007}.
\end{enumerate}
\end{punto}

\begin{notation}
With $\mathcal{L}(M)$ ($\mathcal{L}^{\mathrm{f.i.}}(M)$) we denote the
lattice of (fully invariant) $R$-submodules of $M.$ Moreover, for every $%
L\leq _{R}M$ we set
\begin{equation*}
\mathcal{K}^{_{\mathrm{f.i.}}}(L):=\{\widetilde{L}\leq _{R}M\mid \widetilde{L%
}\subseteq L\text{ and }\widetilde{L}\leq _{R}^{\mathrm{f.i.}}M\}.
\end{equation*}
\end{notation}

\begin{lemma}
\label{ww-fi} Let $L\leq _{R}^{\mathrm{f.i.}}M.$ Then $\mathcal{L}^{_{%
\mathrm{f.i.}}}(L)\subseteq \mathcal{K}^{_{\mathrm{f.i.}}}(L)$ with equality
in case $_{R}M$ is self-injective. In particular, if $_{R}M$ is
self-injective \emph{(}and duo\emph{)}, then $_{R}L$ is self-injective \emph{%
(}and duo\emph{)}.
\end{lemma}

\begin{punto}
By $\mathcal{S}(M)$ ($\mathcal{S}_{\mathrm{f.i.}}(M)$) we denote the
(possibly empty) class of simple $R$-submodules of $M$ (simple $(R,S)$%
-subbimodules of $_{R}M_{S}$, i.e. non-zero fully invariant $R$-submodules
of $M$ that have no non-zero proper fully invariant $R$-submodules). For
every $L\leq _{R}M,$ we set
\begin{equation*}
\begin{tabular}{lll}
$\mathcal{S}(L)$ & $:=$ & $\{K\in \mathcal{S}(M)\mid K\subseteq L\};$ \\
$\mathcal{S}_{\mathrm{f.i.}}(L)$ & $:=$ & $\{K\in \mathcal{S}_{\mathrm{f.i.}%
}(M)\mid K\subseteq L\}.$%
\end{tabular}%
\end{equation*}
\end{punto}

\begin{punto}
Let $L\leq _{R}M.$ We say that $L$ is \emph{essential} or \emph{large} in $%
M, $ and write $L\trianglelefteq M,$ iff $L\cap \widetilde{L}\neq 0$ for
every $0\neq \widetilde{L}\leq M.$ On the other hand, we say $L$ is \emph{%
superfluous} or \emph{small} in $M$ and we write $L\ll M,$ iff $L+\widetilde{%
L}\neq M$ for every $\widetilde{L}\lvertneqq _{R}M.$ With $\mathrm{Max}(R)$
we denote the spectrum of maximal sided ideals of $R$ and with $\mathrm{Spec}%
(R)$ \emph{prime spectrum} of $R;$ moreover, the \emph{prime radical} of $R$
is defined as%
\begin{equation*}
\mathrm{Prad}(R):=\dbigcap\limits_{P\in \mathrm{Spec}(R)}P.
\end{equation*}%
The \emph{socle} of $M$ is defined as
\begin{equation*}
\mathrm{Soc}(M):=\dsum\limits_{L\in \mathcal{S}(M)}L=\dbigcap\limits_{L%
\trianglelefteq M}L\text{ }\;\text{(}:=0\text{ iff }\mathcal{S}%
(M)=\varnothing \text{)}
\end{equation*}%
Call a semisimple module $M$ \emph{completely inhomogenous} iff $M$ is a
direct sum of pairwise non-isomorphic simple submodules.
\end{punto}

\begin{definition}
We say $_{R}M$ is

\emph{colocal} (or \textit{cocyclic} \cite{Wis1991}, \textit{subdirectly
irreducible} \cite{AF1974}) iff $M$ contains a \emph{smallest} non-zero $R$%
-submodule that is contained in every non-zero $R$-submodule of $M$,
equivalently iff $\dbigcap\limits_{0\neq L\leq _{R}M}L\neq 0$;

\emph{uniform}\textbf{,} iff for any $0\neq L_{1},L_{2}\leq _{R}M,$ also $%
L_{1}\cap L_{2}\neq 0$, equivalently iff every non-zero $R$-submodule of $M$
is essential;

\emph{atomic} iff every $0\neq L\leq _{R}M$ contains a simple $R$-submodule,
equivalently iff $\mathrm{Soc}(L)\neq 0$ for every $0\neq L\leq _{R}M$;

\emph{f.i.-atomic} iff $\mathcal{S}_{\mathrm{f.i.}}(L)\neq \varnothing $ for
every $0\neq L\leq _{R}^{\mathrm{f.i.}}M$, equivalently iff $_{R}M_{S}$ is
atomic.
\end{definition}

\begin{exs}
(\cite{HS2010}) Cofinitely generated modules (modules with finitely
generated essential socles; called also \textit{finitely related} modules),
semisimple modules and Artinian modules are atomic. All left modules over
right perfect (e.g. left Artinian) rings are atomic. The Abelian group $%
\mathbb{Z}$ is uniform but not atomic.
\end{exs}

\begin{definition}
We call an $R$-module $M$ an \textit{S-IAD-module} iff $_{R}M$ is
self-injective, atomic and duo.
\end{definition}

\begin{exs}
\label{ex-siad} The following are classes of S-IAD-modules:

\begin{enumerate}
\item self-injective Artinian uniserial modules: Artinian modules are atomic
by \cite{HS2010} and Artinian uniserial modules are duo by \cite{OHS2006};

\item finitely generated self-injective self-cogenerator modules with
commutative endomorphism rings: such modules are duo and finitely
cogenerated by \cite[48.16]{Wis1991}, whence atomic.

\item self-injective duo left modules over right perfect rings;

\item self-injective duo modules that are cofinitely generated (resp.
semisimple, Artinian).
\end{enumerate}
\end{exs}

\subsection*{Topological Spaces}

\qquad In what follows, we fix some definitions and notions for topological
spaces. For further information, the reader might consult any book in
General Topology (e.g. \cite{Bou1966}).

\begin{definition}
We call a topological space $X$ (\textit{countably}) \textit{compact}, iff
every open cover of $X$ has a finite subcover. Countably compact spaces are
also called \textit{Lindelof spaces}. Note that some authors (e.g. \cite%
{Bou1966, Bou1998}) assume that compact spaces are in addition Hausdorff.
\end{definition}

\begin{punto}
We say a topological space ${\mathbf{X}}$ is \textit{Noetherian} (\textit{%
Artinian}), iff every ascending (descending) chain of open sets is
stationary, equivalently iff every descending (ascending) chain of closed
sets is stationary.
\end{punto}

\begin{definition}
(e.g. \cite{Bou1966}, \cite{Bou1998}) A non-empty topological space $\mathbf{%
X}$ is said to be

\begin{enumerate}
\item \emph{ultraconnected}, iff the intersection of any two non-empty
closed subsets is non-empty.

\item \emph{irreducible} (or \emph{hyperconnected}), iff $\mathbf{X}$ is not
the union of two proper \textit{closed} subsets; equivalently, iff the
intersection of any two non-empty open subsets is non-empty.

\item \emph{connected}, iff $\mathbf{X}$ is not the \textit{disjoint} union
of two proper \textit{closed} subsets; equivalently, iff the only subsets of
$\mathbf{X}$ that are \textit{clopen} (i.e. closed and open) are $%
\varnothing $ and $\mathbf{X}$.
\end{enumerate}
\end{definition}

\begin{punto}
(\cite{Bou1966}, \cite{Bou1998}) Let $\mathbf{X}$ be a non-empty topological
space. A non-empty subset ${\mathcal{A}}\subseteq {\mathbf{X}}$ is an
\textit{irreducible set} in $\mathbf{X}$ iff it's an irreducible space
w.r.t. the relative (subspace) topology; in fact, ${\mathcal{A}}\subseteq {%
\mathbf{X}}$ is irreducible iff for any proper closed subsets $A_{1},A_{2}$
of $\mathbf{X}$ we have
\begin{equation*}
{\mathcal{A}}\subseteq A_{1}\cup A_{2}\Rightarrow {\mathcal{A}}\subseteq
A_{1}\text{ or }{\mathcal{A}}\subseteq A_{2}.
\end{equation*}%
A maximal irreducible subspace of $\mathbf{X}$ is called an \emph{%
irreducible component}. An irreducible component of a topological space is
necessarily closed. The irreducible components of a Hausdorff space are just
the singleton sets.
\end{punto}

\begin{definition}
Let ${\mathbf{X}}$ be a topological space and $\mathcal{Y}\subseteq {\mathbf{%
X}}$ be an closed set. A point $y\in {\mathcal{Y}}$ is said to be a \textit{%
generic point}, iff ${\mathcal{Y}}=\overline{\{y\}}$. If every irreducible
closed subset of ${\mathbf{X}}$ has a unique generic point, then we call ${%
\mathbf{X}}$ a \textit{Sober} space.
\end{definition}

\begin{definition}
A collection $\mathcal{G}$ of subsets of a topological space $\mathbf{X}$ is
\emph{locally finite}, iff every point of $\mathbf{X}$ has a neighborhood
that intersects only finitely many elements of $\mathcal{G}.$
\end{definition}

\section{Fully Coprime Submodules}

\qquad As before, $M$ is a non-zero unital left $R$-module. In this section,
we topologize the spectrum $\mathrm{Spec}^{\mathrm{fc}}(M)$ of submodules
that are \emph{fully coprime in }$M.$ For more information on fully
coprimeness, the interested reader can consult \cite{Abu2006} (see also \cite%
{Wij2006} and \cite{WW2009}).

\begin{notation}
For subsets $L\subseteq M$ and $I\subseteq S,$ we set
\begin{equation}
\mathrm{An}(L):=\{f\in S\mid f(L)=0\}\text{ and }\mathrm{Ke}%
(I)=\bigcap_{f\in I}\mathrm{Ker}(f).  \label{AnKe}
\end{equation}
\end{notation}

\begin{punto}
For $R$-submodules $X,Y\leq _{R}M$ we set
\begin{equation*}
X\odot_{M}Y:=\dbigcap \{f^{-1}(Y)\mid f\in \mathrm{An}(X)\}=\dbigcap%
\limits_{f\in \mathrm{An}(X)}\{\mathrm{Ker}(\pi _{Y}\circ f:M\rightarrow
M/Y)\}.
\end{equation*}
If $X\leq _{R}M$ is fully invariant, then $X\odot_{M}Y\leq _{R}M$ is also
fully invariant; and if $Y\leq _{R}M$ is fully invariant, then $X+Y\subseteq
X\odot_{M}Y.$
\end{punto}

\begin{lemma}
\label{inn-ideal}\emph{(See \cite[Lemma 4.9]{Abu2006})} Let $X,Y\leq _{R}^{%
\mathrm{f.i.}}M.$ If $_{R}M$ is self-cogenerator, then
\begin{equation}
X\odot_{M}Y=\mathrm{Ke}(\mathrm{An}(X)\circ ^{op}\mathrm{An}(Y)).
\end{equation}
\end{lemma}

\begin{definition}
We call a non-zero \emph{fully invariant} $R$-submodule $0\neq K\leq _{R}^{%
\mathrm{f.i.}}M$ \emph{fully coprime in}\textbf{\ }$M$ iff for any fully
invariant $R$-submodules $X,Y\leq _{R}^{\mathrm{f.i.}}M:$%
\begin{equation*}
K\leq _{R}X\odot_{M}Y\Rightarrow K\leq _{R}X\text{ or }K\leq _{R}Y.
\end{equation*}
In particular, we say $_{R}M$ is a \emph{fully coprime module} iff for any
fully invariant $R$-submodules $X,Y\leq _{R}^{\mathrm{f.i.}}M:$%
\begin{equation*}
M=X\odot_{M}Y\Rightarrow M=X\text{ or }M=Y.
\end{equation*}
\end{definition}

\begin{proposition}
\emph{(\cite[1.7.3]{Wij2006}, \cite[3.7]{WW2009})} The following are
equivalent:

\begin{enumerate}
\item $_{R}M$ is fully coprime;

\item $M$ is $M/K$-generated for any $K\lvertneqq _{R}^{\mathrm{f.i.}}M;$

\item Any $M$-generated $R$-module is $M/K$-generated for any $K\lvertneqq
_{R}^{\mathrm{f.i.}}M.$
\end{enumerate}
\end{proposition}

\begin{remark}
The definition of fully coprime modules we adopt is a modification of the
definition of prime modules introduced by Bican et. al. \cite{BJKN80}, where
arbitrary submodules are replaced by fully invariant ones (we call such
modules B-\emph{coprime}). In fact, $_{R}M$ is B-coprime if and only if $M$
is generated by each of its non-zero factor modules. Clearly, every
B-coprime module is fully coprime. A duo module is B-coprime if and only if
it is fully coprime. For more details on fully coprime modules, the reader
is referred to \cite{Abu2006} (see also \cite{Wij2006} and \cite{WW2009}).
\end{remark}

\begin{ex}
The Abelian group $\mathbb{Q}$ is fully coprime since it has no non-trivial
fully invariant subgroups. Notice that $\mathbb{_{\mathbb{Z}}Q}$ is not
B-coprime since $\mathbb{Q}$ is not generated by $\mathbb{Q}/\mathbb{Z}$.
\end{ex}

\begin{punto}
We define%
\begin{equation*}
\begin{tabular}{lll}
$\mathrm{Spec}^{\mathrm{fc}}(M)$ & $:=$ & $\{0\neq K\leq _{R}^{\mathrm{f.i.}%
}M\mid K$ is fully coprime in $M\}.$%
\end{tabular}%
\end{equation*}%
We say $_{R}M$ is $\mathrm{fc}$\emph{-coprimeless,} iff $\mathrm{Spec}^{%
\mathrm{fc}}(M)=\varnothing .$ Moreover, for every $R$-submodule $L\leq _{R}M
$ we set
\begin{equation*}
\mathcal{V}^{\mathrm{fc}}(L):=\{K\in \mathrm{Spec}^{\mathrm{fc}}(M)\mid
K\subseteq L\},\ \mathcal{X}^{\mathrm{fc}}(L):=\{K\in \mathrm{Spec}^{\mathrm{%
fc}}(M)\mid K\nsubseteqq L\}
\end{equation*}%
and
\begin{equation*}
\mathrm{Corad}_{M}^{\mathrm{fc}}(L):=\sum\limits_{K\in \mathcal{V}^{\mathrm{%
fc}}(L)}K\text{ \ \ (}:=0,\text{ if }\mathcal{V}^{\mathrm{fc}%
}(L)=\varnothing \text{).}
\end{equation*}%
We say $L\leq _{R}^{\mathrm{f.i.}}M$ is $\mathrm{fc}$\emph{-coradical} iff $%
\mathrm{Corad}_{M}^{\mathrm{fc}}(L)=L.$ In particular, we call $_{R}M$ an $%
\mathrm{fc}$\emph{-coradical module} iff $\mathrm{Corad}_{M}^{\mathrm{fc}%
}(M)=M.$
\end{punto}

\begin{remark}
\label{corad}For any $L_{1}\leq _{R}L_{2}\leq _{R}M$ we have $\mathrm{Corad}%
_{M}^{\mathrm{fc}}(L_{1})\subseteq \mathrm{Corad}_{M}^{\mathrm{fc}}(L_{1}).$
Moreover, for any $L\leq _{R}M$ we have%
\begin{equation*}
\mathrm{Corad}_{M}^{\mathrm{fc}}(\mathrm{Corad}_{M}^{\mathrm{fc}}(L))=%
\mathrm{Corad}_{M}^{\mathrm{fc}}(L).
\end{equation*}
\end{remark}

\begin{punto}
A non-zero fully invariant $R$-submodule $L\leq _{R}^{\mathrm{f.i.}}M$ will
be called $\mathrm{E}$\emph{-prime} iff the ideal $\mathrm{An}(K) \leq S$ is
prime. With $\mathrm{EP}(M)$ we denote the class of $\mathrm{E}$-prime $R$%
-submodules of $M.$
\end{punto}

Recall that $_{R}M$ is said to be \emph{intrinsically injective} iff $%
\mathrm{AnKe}(I)=I$ for every finitely generated right ideal $I\leq S$.
Every self-injective $R$-module is intrinsically injective (e.g. \cite[28.1]%
{Wis1991}).

\begin{proposition}
\label{corad=}\emph{(\cite[Proposition 4.12]{Abu2006})} If $_{R}M$ is
self-cogenerator, then $\mathrm{EP}(M)\subseteq \mathrm{Spec}^{\mathrm{fc}%
}(M)$ with equality if $_{R}M$ is intrinsically injective. If moreover $S$
is right Noetherian, then
\begin{equation*}
\mathrm{Prad}(S)=\mathrm{An}(\mathrm{Corad}^{\mathrm{fc}}(M))\text{ and }%
\mathrm{Corad}^{\mathrm{fc}}(M)=\mathrm{Ke}(\mathrm{Prad}(S)).
\end{equation*}
\end{proposition}

\begin{proposition}
\label{ro-inn}\emph{(\cite[Proposition 4.7]{Abu2006})}\ Let $0\neq L\leq
_{R}^{\mathrm{f.i.}}M.$ Then
\begin{equation*}
\mathcal{L}^{\mathrm{f.i.}}(L)\cap \mathrm{Spec}^{\mathrm{fc}}(M)\subseteq
\mathrm{Spec}^{\mathrm{fc}}(L),
\end{equation*}
with equality if $_{R}M$ is self-injective. In particular, if $_{R}M$ is
self injective then
\begin{eqnarray*}
L\text{ is fully coprime in }M\Leftrightarrow \text{ }_{R}L\text{ is a fully
coprime module.}
\end{eqnarray*}
\end{proposition}

\begin{remark}
\label{sub-cop}Every $L\in \mathcal{S}_{\mathrm{f.i.}}(M)$ is trivially a
fully coprime $R$-module. If $_{R}M$ is self-injective, then $\mathcal{S}_{%
\mathrm{f.i.}}(M)\subseteq \mathrm{Spec}^{\mathrm{fc}}(M)$ by Proposition %
\ref{ro-inn}; hence if, in addition, $M$ is f.i.-atomic, then for every $%
L\leq _{R}^{\mathrm{f.i.}}M$ we have $\varnothing \neq \mathcal{S}_{\mathrm{%
f.i.}}(L)\subseteq \mathrm{Spec}^{\mathrm{fc}}(L)=\mathcal{V}^{\mathrm{fc}%
}(L)\subseteq \mathrm{Spec}^{\mathrm{fc}}(M).$
\end{remark}

\begin{definition}
Let $K\in \mathrm{Spec}^{\mathrm{fc}}(M).$ We say that $K$ is \emph{maximal
under }$L,$ where $0\neq L\leq _{R}^{\mathrm{f.i.}}M$ iff $K$ is a maximal
element of $\mathcal{V}^{\mathrm{fc}}(L)$, equivalently $K\subseteq L$ and
there is no $\widetilde{K}\in \mathrm{Spec}^{\mathrm{fc}}(M)$ that is
contained in $L$ and contains $K$ strictly. We say that $K$ is \emph{maximal
in }$\mathrm{Spec}^{\mathrm{fp}}(M)$ iff $K$ is maximal under $M.$
\end{definition}

\begin{lemma}
\label{fc-max}Let $M$ be self-injective and f.i.-atomic. For every $0\neq
L\leq _{R}^{\mathrm{f.i.}}M$ there exists $K\in \mathrm{Spec}^{\mathrm{fc}%
}(M)$ which is maximal under $L.$ In particular, $\mathrm{Spec}^{\mathrm{fc}%
}(M)$ has a maximal element.
\end{lemma}

\begin{Beweis}
Let $0\neq L\leq _{R}^{\mathrm{f.i.}}M.$ By Remark \ref{sub-cop}, $%
\varnothing \neq \mathcal{S}_{\mathrm{f.i.}}(L)\subseteq \mathcal{V}^{%
\mathrm{fc}}(L)$. Let
\begin{equation*}
K_{1}\subseteq K_{2}\subseteq \cdots \subseteq K_{n}\subseteq
K_{n+1}\subseteq \cdots
\end{equation*}%
be an ascending chain in $\mathcal{V}^{\mathrm{fc}}(L)$ and set $\widetilde{K%
}:=\bigcup\limits_{i=1}^{\infty }K_{i}.$ Suppose that there exist two fully
invariant $R$-submodules $L_{1},L_{2}\leq _{R}^{\mathrm{f.i.}}M$ with $%
\widetilde{K}\subseteq L_{1}\odot _{M}L_{2}$ but $\widetilde{K}\nsubseteqq
L_{1}$ and $\widetilde{K}\nsubseteqq L_{2}.$ Then $K_{n_{1}}\nsubseteqq L_{1}
$ for some $n_{1}$ and $K_{n_{2}}\nsubseteqq L_{2}$ for some $n_{2}.$
Setting $n:=\max \{n_{1},n_{2}\},$ we have $K_{n}\subseteq L_{1}\odot
_{M}L_{2}$ while $K_{n}\nsubseteqq L_{1}$ and $K_{n}\nsubseteqq L_{2}$, a
contradiction. So, $\widetilde{K}\in \mathcal{V}^{\mathrm{fc}}(L).$ By
Zorn's Lemma, $\mathcal{V}^{\mathrm{fc}}(L)$ has a maximal element. In
particular, $\mathrm{Spec}^{\mathrm{fc}}(M)=\mathcal{V}^{\mathrm{fc}}(M)$
has a maximal element.$\blacksquare $
\end{Beweis}

\qquad We introduce now a notion for modules which will prove to be useful
in the sequel. Moreover, this notion seems to be of independent interest
(see the survey in \cite{Smi}):

\begin{punto}
We say that $_{R}M$ has the \emph{min-property} iff for any simple $R$%
-submodule $L\in \mathcal{S}(M)$ we have $L\nsubseteqq L_{e},$ where
\begin{equation}
L_{e}:=\dsum\limits_{K\in \mathcal{S}(M)\backslash \{L\}}K\text{ \ \ (}:=0%
\text{, if }\mathcal{S}(M)=\{L\}\text{).}  \label{A_L}
\end{equation}%
Since simple modules are cyclic, $_{R}M$ has the min-property if and only if
for any $L\in \mathcal{S}(M)$ and any \textit{finite} subset $\{L_{1},\cdots
,L_{n}\}\subseteq \mathcal{S}(M)\setminus \{L\}$, we have $L\nsubseteqq
\sum_{i=1}^{n}L_{i}$. It is obvious that $_{R}M$ has the min-property if and
only if the class $\mathcal{S}(M)$ of simple $R$-submodules is \textit{%
independent} in the sense of \cite[page 8]{CLVW2006}. By \cite[Theorem 2.3]%
{Smi}, $_{R}M$ has the min-property if and only if all distinct simple $R$%
-submodules of $M$ are non-isomorphic (i.e. $\mathrm{Soc}(M)$ is completely
inhomogenous).
\end{punto}

\begin{exs}
\begin{enumerate}
\item Every $R$-module with at most one simple $R$-submodule (e.g. a colocal
$R$-module) has the min-property.

\item Let $R$ have the property that $R/P$ is left Artinian for every
primitive left ideal $P$ of $R$ (e.g. $R$ is a commutative ring, or a
PI-ring, or a semilocal ring). Let $\{U_{\lambda }\}_{\lambda \in \Lambda }$
be any collection of pairwise non-isomorphic simple $R$-modules and consider
$M:=\dprod\limits_{\lambda \in \Lambda }U_{\lambda }.$ Then $\mathrm{Soc}%
(_{R}M)=\dbigoplus\limits_{\lambda \in \Lambda }U_{\lambda }$ is completely
inhomogenous and so $_{R}M$ has the min-property \cite[Proposition 2.6]{Smi}.

\item Any finitely generated Artinian module over a commutative ring has the
min-property if and only if it is a comultiplication module \cite[Theorem
3.11]{AS}.
\end{enumerate}
\end{exs}

\begin{lemma}
\label{min} If $_{R}M$ is self-injective and duo, then $_{R}M$ has the
min-property.
\end{lemma}

\begin{Beweis}
Assume, without loss of generality, that $\left\vert \mathcal{S}%
(M)\right\vert \geq 2.$ Let $L\leq _{R}M$ be simple and suppose that $%
L\subseteq \sum_{i=1}^{n}L_{i}$ for $\{L_{1},\cdots ,L_{n}\}\subseteq
\mathcal{S}(M)\setminus \{L\}$. Then $L\subseteq L_{1}\odot _{M}\cdots \odot
_{M}L_{n}.$ Since every simple $R$-submodule of $M$ is fully coprime in $M$
(see Remark \ref{sub-cop}), induction yields that $L=L_{i}$ for some $%
i=1,\cdots ,n$, a contradiction.$\blacksquare $
\end{Beweis}

\subsection*{Top$^{\mathrm{fc}}$-modules}

\begin{notation}
Set
\begin{equation*}
\begin{tabular}{lllllll}
$\xi ^{\mathrm{fc}}(M)$ & $:=$ & $\{\mathcal{V}^{\mathrm{fc}}(L)\mid L\leq_R
M\};$ &  & $\xi _{f.i.}^{\mathrm{fc}}(M)$ & $:=$ & $\{\mathcal{V}^{\mathrm{fc%
}}(L)\mid L\leq _{R}^{\mathrm{f.i.}}M\};$ \\
$\tau ^{\mathrm{fc}}(M)$ & $:=$ & $\{\mathcal{X}^{\mathrm{fc}}(L)\mid L\leq
_{R}M\};$ &  & $\tau _{f.i.}^{\mathrm{fc}}(M)$ & $:=$ & $\{\mathcal{X}^{%
\mathrm{fc}}(L)\mid L\leq _{R}^{\mathrm{f.i.}}M\};$ \\
$\mathbf{Z}^{\mathrm{fc}}(M)$ & $:=$ & $(\mathrm{Spec}^{\mathrm{fc}}(M),\tau
^{\mathrm{fc}}(M));$ &  & $\mathbf{Z}_{f.i.}^{\mathrm{fc}}(M)$ & $:=$ & $(%
\mathrm{Spec}^{\mathrm{fc}}(M),\tau _{f.i.}^{\mathrm{fc}}(m)).$%
\end{tabular}%
\end{equation*}
\end{notation}

\begin{lemma}
\label{Properties}

\begin{enumerate}
\item $\mathcal{V}^{\mathrm{fc}}(0)=\emptyset $ and $\mathcal{V}^{\mathrm{fc}%
}(M)=\mathrm{Spec}^{\mathrm{fc}}(M).$

\item If $\{L_{\lambda }\}_{\Lambda }\subseteq \mathcal{L}(M),$ then
\begin{equation*}
\bigcap\limits_{\Lambda }\mathcal{V}^{\mathrm{fc}}(L_{\lambda })=\mathcal{V}%
^{\mathrm{fc}}(\bigcap\limits_{\Lambda }L_{\lambda }).
\end{equation*}

\item If $L,\widetilde{L}\in \mathcal{L}^{\mathrm{f.i.}}(M),$ then
\begin{equation*}
\mathcal{V}^{\mathrm{fc}}(L)\cup \mathcal{X}^{\mathrm{fc}}(\widetilde{L})=%
\mathcal{V}^{\mathrm{fc}}(L+\widetilde{L})=\mathcal{V}^{\mathrm{fc}%
}(L\odot_{M}\widetilde{L}).
\end{equation*}
\end{enumerate}
\end{lemma}

\begin{Beweis}
Statements \textquotedblleft 1\textquotedblright\ and \textquotedblleft
2\textquotedblright\ and the inclusions $\mathcal{V}^{\mathrm{fc}}(L)\cup
\mathcal{V}^{\mathrm{fc}}(\widetilde{L})\subseteq \mathcal{V}^{\mathrm{fc}%
}(L+\widetilde{L})\subseteq \mathcal{V}^{\mathrm{fc}}(L\odot _{M}\widetilde{L%
})$ in (3) are obvious. Let $K\in \mathcal{V}^{\mathrm{fc}}(L\odot _{M}%
\widetilde{L}),$ so that $K\subseteq L\odot _{M}\widetilde{L}.$ Since $K$ is
fully coprime in $M,$ we have $K\subseteq L$ whence $K\in \mathcal{V}^{%
\mathrm{fc}}(L),$ or $K\subseteq \widetilde{L}$ whence $K\in \mathcal{V}^{%
\mathrm{fc}}(\widetilde{L}).\blacksquare $
\end{Beweis}

For an arbitrary $R$-module $M,$ the set $\xi ^{\mathrm{fc}}(M)$ is not
necessarily closed under finite unions. This motivates the following

\begin{definition}
We call $_{R}M$ an \emph{top}$^{\mathrm{fc}}$\emph{-module} iff $\xi ^{%
\mathrm{fc}}(M)$ is closed under finite unions.
\end{definition}

\begin{notation}
For any $\mathcal{A}\subseteq \mathrm{Spec}^{\mathrm{fc}}(M)$ set
\begin{equation*}
\mathcal{H}(\mathcal{A}):=\dsum_{K\in \mathcal{A}}K\text{ }\;\;\text{(}:=0,%
\text{ if }\mathcal{A}=\varnothing \text{).}
\end{equation*}
\end{notation}

\qquad As a direct consequence of Lemma \ref{Properties} we get

\begin{theorem}
\label{Topology}$\mathbf{Z}_{f.i.}^{\mathrm{fc}}(M):=(\mathrm{Spec}^{\mathrm{%
fc}}(M),\tau _{f.i.}^{\mathrm{fc}})$ is a topological space. In particular,
if $_R M$ is duo, then $M$ is an top$^{\mathrm{fc}}$-module and $\mathbf{Z}^{%
\mathrm{fc}}(M):=(\mathrm{Spec}^{\mathrm{fc}}(M),\tau ^{\mathrm{fc}}(M))$ is
a topological space.
\end{theorem}

\begin{lemma}
\label{closure} Let $_{R}M$ be a top$^{\mathrm{fc}}$-module. The closure of
any subset $\mathcal{A}\subseteq \mathrm{Spec}^{\mathrm{fc}}(M)$ is
\begin{equation}
\overline{\mathcal{A}}=\mathcal{V}^{\mathrm{fc}}(\mathcal{H}(\mathcal{A})%
\mathcal{)}.  \label{fc-closure}
\end{equation}
\end{lemma}

\begin{Beweis}
Let $\mathcal{A}\subseteq \mathrm{Spec}^{\mathrm{fc}}(M).$ Since $\mathcal{A}%
\subseteq \mathcal{V}^{\mathrm{fc}}(\mathcal{H}(\mathcal{A}))$ and $\mathcal{%
V}^{\mathrm{fc}}(\mathcal{H}(\mathcal{A}))$ is a closed set, we have $%
\overline{\mathcal{A}}\subseteq \mathcal{V}^{\mathrm{fc}}(\mathcal{H}(%
\mathcal{A})).$ On the other hand, suppose that $H\in \mathcal{V}^{\mathrm{fc%
}}(\mathcal{H}(\mathcal{A}))\backslash \mathcal{A}$ and let $\mathcal{X}^{%
\mathrm{fc}}(L)$ be a neighborhood of $H,$ so that $H\nsubseteqq L.$ Then
there exists $W\in \mathcal{A}$ with $W\nsubseteqq L$ (otherwise $H\subseteq
\mathcal{H}(\mathcal{A})\subseteq L,$ a contradiction), i.e. $W\in \mathcal{X%
}^{\mathrm{fc}}(L)\cap (\mathcal{A}\backslash \{H\})$ is a cluster point of $%
\mathcal{A}$ Consequently, $\overline{\mathcal{A}}=\mathcal{V}^{\mathrm{fc}}(%
\mathcal{H}(\mathcal{A})).\blacksquare $
\end{Beweis}

\begin{remarks}
\label{simple-char}Let $M$ be a top$^{\mathrm{fc}}$-module and consider the
Zariski topology $\mathbf{Z}^{\mathrm{fc}}(M):=(\mathrm{Spec}^{\mathrm{fc}%
}(M),\tau ^{\mathrm{fc}}(M)).$

\begin{enumerate}
\item $\mathbf{Z}^{\mathrm{fc}}(M)$ is a $T_{0}$ (Kolmogorov) space.

\item If $_{R}M$ is duo, then $\mathcal{B}:=\{\mathcal{X}^{\mathrm{fc}%
}(L)\mid L\leq _{R}M$ is finitely generated$\}$ is a basis of open sets for $%
\mathbf{Z}^{\mathrm{fc}}(M):$ any $K\in \mathrm{Spec}^{\mathrm{fc}}(M)$ is
contained in some $\mathcal{X}^{\mathrm{fc}}(L)$ for some finitely generated
$R$-submodule $L\leq _{R}M$ (e.g. $L=0$). Moreover, if $L_{1},L_{2}\leq _{R}M
$ are finitely generated and $K\in \mathcal{X}^{\mathrm{fc}}(L_{1})\cap
\mathcal{X}^{\mathrm{fc}}(L_{2}),$ then $L:=L_{1}+L_{2}\leq _{R}M$ is also
finitely generated and we have $K\in \mathcal{X}^{\mathrm{fc}}(L)=\mathcal{X}%
^{\mathrm{fc}}(L_{1})\cap \mathcal{X}^{\mathrm{fc}}(L_{2}).$

\item If $L\in \mathrm{Spec}^{\mathrm{fc}}(M),$ then $\overline{\{L\}}=%
\mathcal{V}^{\mathrm{fc}}(L).$ In particular, for any $K\in \mathrm{Spec}^{%
\mathrm{fc}}(M):$%
\begin{equation*}
K\in \overline{\{L\}}\Leftrightarrow K\subseteq L.
\end{equation*}

\item If $_{R}M$ is self-injective and duo, then
\begin{equation}
{\mathcal{X}}^{\mathrm{fc}}(L)=\emptyset \Rightarrow \mathrm{Soc}%
(M)\subseteq L.  \label{Soc}
\end{equation}%
The converse of (\ref{Soc}) holds if, for example, $\mathcal{S}(M)=\mathrm{%
Spec}^{\mathrm{fc}}(M).$

\item Let $_{R}M$ be a S-IAD-module. For every $L\leq _{R}M$ we have ${%
\mathcal{V}}^{\mathrm{fc}}(L)=\emptyset $ if and only if $L=0.$

\item Let $_{R}M$ be self-injective and $0\neq L\leq _{R}^{\mathrm{f.i.}}M.$
The embedding
\begin{equation*}
\mathbf{\iota }:\mathrm{Spec}^{\mathrm{fc}}(L)\rightarrow \mathrm{Spec}^{%
\mathrm{fc}}(M)
\end{equation*}%
is continuous: this follows from Proposition \ref{ro-inn}, which implies
that that $\iota ^{-1}(\mathcal{V}^{\mathrm{fc}}(N))=\mathcal{V}^{\mathrm{fc}%
}(N\cap L)$ for every $R$-submodule $N\leq _{R}M.$

\item Let $M\overset{\theta }{\simeq }N$ be an isomorphism of non-zero $R$%
-modules. Then we have a bijection
\begin{equation*}
\mathrm{Spec}^{\mathrm{fc}}(M)\longleftrightarrow \mathrm{Spec}^{\mathrm{fc}%
}(N).
\end{equation*}%
In particular, we have $\theta (\mathrm{Corad}_{M}^{\mathrm{fc}}(M))=\mathrm{%
Corad}_{N}^{\mathrm{fc}}(N).$ Moreover, we have a \emph{homeomorphism }$%
\mathbf{Z}^{\mathrm{fc}}(M)\approx \mathbf{Z}^{\mathrm{fc}}(N).$
\end{enumerate}
\end{remarks}

\begin{notation}
Set
\begin{equation*}
\mathbf{CL}(\mathbf{Z}^{\mathrm{fc}}(M)):=\{\mathcal{A}\subseteq \mathrm{Spec%
}^{\mathrm{fc}}(M)\mid \mathcal{A}=\overline{\mathcal{A}}\}\text{ and }%
\mathcal{CR}^{\mathrm{fc}}(M):=\{L\leq _{R}M\mid \mathrm{Corad}_{M}^{\mathrm{%
fc}}(L)=L\}.
\end{equation*}
\end{notation}

\begin{theorem}
\label{11} Let $M$ be a top$^{\mathrm{fc}}$-module.

\begin{enumerate}
\item We have an order-preserving bijection
\begin{equation}  \label{bij}
\mathcal{CR}^{\mathrm{fc}}(M){\longleftrightarrow }\mathbf{CL}(\mathbf{Z}^{%
\mathrm{fc}}(M)),\text{ }L\mapsto \mathcal{V}^{\mathrm{fc}}(L).
\end{equation}

\item $\mathbf{Z}^{\mathrm{fc}}(M)$ is Noetherian if and only if $_R M$
satisfies the DCC condition on fc-coradical submodules.

\item $\mathbf{Z}^{\mathrm{fc}}(M)$ is Artinian if and only if $_R M$
satisfies the ACC condition on fc-coradical submodules.
\end{enumerate}
\end{theorem}

\begin{Beweis}
Notice that for any $L\leq _{R}M$ we have $\mathcal{V}^{\mathrm{fc}}(L)\in
\mathbf{CL}(\mathbf{Z}^{\mathrm{fc}}(M))$ and $\mathrm{Corad}_{M}^{\mathrm{fc%
}}(L)\in \mathcal{CR}^{\mathrm{fc}}(M)$ by Remark \ref{corad}. Consider now
\begin{equation*}
\psi :\mathbf{CL}(\mathbf{Z}^{\mathrm{fp}}(M))\rightarrow \mathcal{CR}^{%
\mathrm{fc}}(M),\text{ }\mathcal{V}^{\mathrm{fc}}(L)\mapsto \mathrm{Corad}%
_{M}^{\mathrm{fc}}(L).
\end{equation*}%
For every $L\in \mathcal{CR}^{\mathrm{fc}}(M)$ we have
\begin{equation*}
\psi (\mathcal{V}^{\mathrm{fc}}(L))=\mathrm{Corad}_{M}^{\mathrm{fc}}(L)=L.
\end{equation*}%
Moreover, for every $\mathcal{A}:=\mathcal{V}^{\mathrm{fc}}(K)\in \mathbf{CL}%
(\mathbf{Z}^{\mathrm{fc}}(M))$ we have
\begin{equation*}
\mathcal{V}^{\mathrm{fc}}(\psi ({\mathcal{A}}))=\mathcal{V}^{\mathrm{fc}%
}(\psi (\mathcal{V}^{\mathrm{fc}}(K)))=\mathcal{V}^{\mathrm{fc}}(\mathrm{%
Corad}_{M}^{\mathrm{fc}}(K))=\mathcal{V}^{\mathrm{fc}}(\mathcal{H}(\mathcal{A%
}))=\overline{\mathcal{A}}=\mathcal{A}.
\end{equation*}%
Notice now that (2) and (3) follow directly from (1) and so we are done.$%
\blacksquare $
\end{Beweis}

\begin{theorem}
\label{noth-art} Let $M$ be a top$^{\mathrm{fc}}$-module. If $_{R}M$ is
Artinian \emph{(}Noetherian\emph{)}, then $\mathbf{Z}^{\mathrm{fc}}(M)$ is
Noetherian \emph{(}Artinian\emph{)}.
\end{theorem}


\begin{proposition}
\label{duo-irr}Let $_{R}M$ be duo. Then ${\mathcal{A}}\subseteq \mathrm{Spec}%
^{\mathrm{fc}}(M)$ is irreducible if and only if ${\mathcal{H}}({\mathcal{A}}%
)$ is fully coprime in $M$.
\end{proposition}

\begin{Beweis}
Let $_R M$ be duo and ${\mathcal{A}} \subseteq \mathrm{Spec}^{\mathrm{fc}%
}(M) $.

($\Rightarrow $)\ Assume that ${\mathcal{A}}$ is irreducible. By definition,
${\mathcal{A}}\neq \emptyset $ and so ${\mathcal{H}}({\mathcal{A}})\neq 0.$
Suppose that ${\mathcal{H}}({\mathcal{A}})$ is not fully coprime in $M$, so
that there exist $R$-submodules $X,Y\leq _{R}M$ with ${\mathcal{H}}({%
\mathcal{A}})\subseteq X\odot _{M}Y$ but ${\mathcal{H}}({\mathcal{A}}%
)\nsubseteqq X$ and ${\mathcal{H}}({\mathcal{A}})\nsubseteqq Y.$ It follows
that ${\mathcal{A}}\subseteq \mathcal{V}^{\mathrm{fc}}(X\odot _{M}Y)=%
\mathcal{V}^{\mathrm{fc}}(X)\cup \mathcal{V}^{\mathrm{fc}}(Y)$ a union of
two proper closed subsets, a contradiction. Consequently, ${\mathcal{H}}({%
\mathcal{A}})$ is fully coprime in $M.$

($\Leftarrow $) Assume that ${\mathcal{H}}({\mathcal{A}})\in \mathrm{Spec}^{%
\mathrm{fc}}(M).$ In particular, ${\mathcal{H}}({\mathcal{A}})\neq 0$ and so
${\mathcal{A}}\neq \varnothing .$ Suppose that ${\mathcal{A}}\subseteq
\mathcal{V}^{\mathrm{fc}}(L_{1})\cup \mathcal{V}^{\mathrm{fc}}(L_{2})=%
\mathcal{V}^{\mathrm{fc}}(L_{1}\odot _{M}L_{2})$ for some $R$-submodules $%
0\neq L_{1},L_{2}\leq _{R}M.$ It follows that ${\mathcal{H}}({\mathcal{A}}%
)\subseteq \mathcal{V}^{\mathrm{fc}}(L_{1}\odot _{M}L_{2})$ and it follows
from our assumption that that ${\mathcal{H}}({\mathcal{A}})\subseteq L_{1}$
so that $\mathcal{A}\subseteq \mathcal{V}^{\mathrm{fc}}(L_{1})$ or ${%
\mathcal{H}}({\mathcal{A}})\subseteq L_{2}$ so that $\mathcal{A}\subseteq
\mathcal{V}^{\mathrm{fc}}(L_{2})$. Consequently, $\mathcal{A}$ is not the
union of two \emph{proper} closed subsets, i.e. $\mathcal{A}$ is irreducible.%
$\blacksquare $
\end{Beweis}

\begin{theorem}
\label{corad-fc} Let $_R M$ be duo.

\begin{enumerate}
\item $\mathrm{Spec}^{\mathrm{fc}}(M)$ is irreducible if and only if $%
\mathrm{Corad}^{\mathrm{fc}}_M(M)$ is fully coprime in $M$.

\item If $_R M$ is self-injective, then ${\mathcal{S}}(M)$ is irreducible if
and only if $\mathrm{Soc}(M)$ is fully coprime in $M$.
\end{enumerate}
\end{theorem}

\begin{ex}
If ${\mathcal{A}}\subset \mathrm{Spec}^{\mathrm{fc}}(M)$ is a chain, then ${%
\mathcal{A}}$ is irreducible. In particular, if $_{R}M$ is uniserial, then%
\textrm{\ }$\mathrm{Spec}^{\mathrm{fc}}(M)$ is irreducible.
\end{ex}

\begin{proposition}
\label{max-irr}Let $_{R}M$ be duo. The bijection \emph{(}\ref{bij}\emph{)}
restricts to bijections:
\begin{equation*}
\mathrm{Spec}^{\mathrm{fc}}(M)\longleftrightarrow \{{\mathcal{A}}\mid {%
\mathcal{A}}\subseteq \mathrm{Spec}^{\mathrm{fc}}(M)\text{ is an irreducible
closed subset}\}
\end{equation*}%
and
\begin{equation*}
\mathrm{Max}(\mathrm{Spec}^{\mathrm{fc}}(M))\longleftrightarrow \{{\mathcal{A%
}}\mid {\mathcal{A}}\subseteq \mathrm{Spec}^{\mathrm{fc}}(M)\text{ is an
irreducible component}\}.
\end{equation*}
\end{proposition}

\begin{Beweis}
Recall the bijection $\mathcal{CR}^{\mathrm{fc}}(M)\overset{\mathcal{V}^{%
\mathrm{fc}}(-)}{\longrightarrow }\mathbf{CL}(\mathbf{Z}^{\mathrm{fc}}(M)).$

Let $K\in \mathrm{Spec}^{\mathrm{fc}}(M).$ Then $K={\mathcal{H}}(\mathcal{V}%
^{\mathrm{fc}}(K))$ and it follows that the closed set $\mathcal{V}^{\mathrm{%
fc}}(K)$ is irreducible by Proposition \ref{duo-irr}. On the other hand, let
${\mathcal{A}}\subseteq \mathrm{Spec}^{\mathrm{fc}}(M)$ be a closed
irreducible subset. Then ${\mathcal{A}}=\mathcal{V}^{\mathrm{fc}}(L)$ for
some $L\leq _{R}M.$ Notice that ${\mathcal{H}}({\mathcal{A}})$ is fully
coprime in $M$ by Proposition \ref{duo-irr} and that ${\mathcal{A}}=%
\overline{{\mathcal{A}}}=\mathcal{V}^{\mathrm{fc}}({\mathcal{H}}({\mathcal{A}%
})).$

On the other hand, notice that $\mathrm{Spec}^{\mathrm{fc}}(M)$ has maximal
elements by Lemma \ref{fc-max}. If $K$ is maximal in $\mathrm{Spec}^{\mathrm{%
fc}}(M),$ then clearly $\mathcal{V}^{\mathrm{fc}}(K)$ is an irreducible
component of $\mathrm{Spec}^{\mathrm{fc}}(M)$ by \textquotedblleft
1\textquotedblright . Conversely, let $\mathcal{Y}$ be an irreducible
component of $\mathrm{Spec}^{\mathrm{fc}}(M).$ Then $\mathcal{Y}$ is closed.
By \textquotedblleft 1\textquotedblright , $\mathcal{Y}=\mathcal{V}^{\mathrm{%
fc}}(L)$ for some $L\in \mathrm{Spec}^{\mathrm{fc}}(M).$ Suppose that $L$ is
not maximal in $\mathrm{Spec}^{\mathrm{fc}}(M)$, so that there exists $K\in
\mathrm{Spec}^{\mathrm{fc}}(M)$ such that $L\subsetneqq K\subseteq M.$ It
follows that $\mathcal{V}^{\mathrm{fc}}(L)\subsetneqq \mathcal{V}^{\mathrm{fc%
}}(L)$, a contradiction since $\mathcal{V}^{\mathrm{fc}}(K)\subseteq \mathrm{%
Spec}^{\mathrm{fc}}(M)$ is irreducible by \textquotedblleft
1\textquotedblright . We conclude that $L$ is maximal in $\mathrm{Spec}^{%
\mathrm{fc}}(M).\blacksquare $
\end{Beweis}

\begin{corollary}
Let $_R M$ be duo. Then $\mathrm{Spec}^{\mathrm{fp}}(M)$ is a Sober space.
\end{corollary}

\begin{Beweis}
Let ${\mathcal{A}} \subseteq \mathrm{Spec}^{\mathrm{fc}}(M)$ be a closed
irreducible subset. By Proposition \ref{max-irr}, ${\mathcal{A}} = {\mathcal{%
V}}^{\mathrm{fc}}(K)$ for some $K \in \mathrm{Spec}^{\mathrm{fp}}(M)$. It
follows that
\begin{equation*}
{\mathcal{A}} = \overline{{\mathcal{A}}} = {\mathcal{V}}^{\mathrm{fc}}({%
\mathcal{H}}({\mathcal{A}})) = {\mathcal{V}}^{\mathrm{fc}}(K) = \overline{%
\{K\}},
\end{equation*}
i.e. $K$ is a generic point for ${\mathcal{A}}$. If $L$ is a generic point
of ${\mathcal{A}}$, then it follows that ${\mathcal{V}}^{\mathrm{fc}}(K) = {%
\mathcal{V}}^{\mathrm{fc}}(L)$ whence $K = L$ since $K, L \in \mathrm{Spec}^{%
\mathrm{fc}}(M)$.$\blacksquare$
\end{Beweis}

\begin{proposition}
\label{uniform} Let $_R M$ be a S-IAD-module. Then $_R M$ is uniform if and
only if $\mathrm{Spec}^{\mathrm{fc}}(M)$ is ultraconnected.
\end{proposition}

\begin{Beweis}
Assume that $_{R}M$ is uniform. If $\mathcal{V}^{\mathrm{fc}}(L_{1}),$ $%
\mathcal{V}^{\mathrm{fc}}(L_{2})\subseteq \mathrm{Spec}^{\mathrm{fc}}(M)$
are any two non-empty closed subsets, then $L_{1}\neq 0\neq L_{2}$ and so $%
\mathcal{V}^{\mathrm{fc}}(L_{1})\cap \mathcal{V}^{\mathrm{fc}}(L_{2})=%
\mathcal{V}^{\mathrm{fc}}(L_{1}\cap L_{2})\neq \varnothing ,$ since $%
L_{1}\cap L_{2}\neq 0$ contains by assumption some simple $R$-submodule
whence which is indeed fully coprime in $M.$ Conversely, assume that the
intersection of any two non-empty closed subsets of $\mathrm{Spec}^{\mathrm{%
fc}}(M)$ is non-empty. Let $0\neq L_{1},L_{2}\leq _{R}M,$ so that $\mathcal{V%
}^{\mathrm{fc}}(L_{1})\neq \varnothing \neq \mathcal{V}^{\mathrm{fc}}(L_{2}).
$ By assumption $\mathcal{V}^{\mathrm{fc}}(L_{1}\cap L_{2})=\mathcal{V}^{%
\mathrm{fc}}(L_{1})\cap \mathcal{V}^{\mathrm{fc}}(L_{2})\neq \varnothing ,$
hence $L_{1}\cap L_{2}\neq 0.$ Consequently, $_{R}M$ is uniform.$%
\blacksquare $
\end{Beweis}

\begin{theorem}
\label{compact} Let $_R M$ be a S-IAD-module.

\begin{enumerate}
\item If $\mathcal{S}(M)$ is countable, then $\mathbf{Z}^{\mathrm{fc}}(M)$
is countably compact.

\item If $\mathcal{S}(M)$ is finite, then $\mathbf{Z}^{\mathrm{fc}}(M)$ is
compact.
\end{enumerate}
\end{theorem}

\begin{Beweis}
We prove only \textquotedblleft 1\textquotedblright , since
\textquotedblleft 2\textquotedblright\ can be proved similarly. Assume that $%
\mathcal{S}(M)=\{N_{\lambda _{k}}\}_{k\geq 1}$ is countable. Let $\{\mathcal{%
X}^{\mathrm{fc}}(L_{\alpha })\}_{\alpha \in I}$ be an open cover of $\mathrm{%
Spec}^{\mathrm{fc}}(M)$, i.e. $\mathrm{Spec}^{\mathrm{fc}}(M)\subseteq
\dbigcup\limits_{\alpha \in I}\mathcal{X}^{\mathrm{fc}}(L_{\alpha })$. Since
$\mathcal{S}(M)\subseteq \mathrm{Spec}^{\mathrm{fc}}(M),$ we can pick for
each $k\geq 1,$ some $\alpha _{k}\in I$ such that $N_{\lambda
_{k}}\nsubseteqq L_{\alpha _{k}}.$ Suppose $\dbigcap\limits_{k\geq
1}L_{\alpha _{k}}\neq 0.$ Since $_{R}M$ is atomic, there exists some simple $%
R$-submodule $0\neq N\subseteq \dbigcap\limits_{k\geq 1}L_{\alpha _{k}},$ a
contradiction since $N=N_{\lambda _{k}}\nsubseteqq L_{\alpha _{k}}$ for some
$k\geq 1$. Hence $\dbigcap\limits_{k\geq 1}L_{\alpha _{k}}=0$ and we
conclude that $\mathrm{Spec}^{\mathrm{fc}}(M)=\mathcal{X}^{\mathrm{fc}%
}(\dbigcap\limits_{k\geq 1}L_{\alpha _{k}})=\dbigcup\limits_{k\geq 1}%
\mathcal{X}^{\mathrm{fc}}(L_{\alpha _{k}})$, i.e. $\{\mathcal{X}^{\mathrm{fc}%
}(L_{\alpha _{k}})\mid k\geq 1\}\subseteq \{\mathcal{X}^{\mathrm{fc}%
}(L_{\alpha })\}_{\alpha \in I}$ is a countable subcover.$\blacksquare $
\end{Beweis}

\begin{proposition}
\label{it-irr} Let $_R M$ be duo and assume that $\mathrm{Spec}^{\mathrm{fc}%
}(M) = {\mathcal{S}}(M).$

\begin{enumerate}
\item If $_R M$ has the min-property, then $\mathrm{Spec}^{\mathrm{fc}}(M)$
is discrete.

\item $M$ has a unique simple $R$-submodule if and only if $_R M$ has the
min-property and $\mathrm{Spec}^{\mathrm{fc}}(M)$ is connected.
\end{enumerate}
\end{proposition}

\begin{Beweis}
\begin{enumerate}
\item If $_{R}M$ has the min-property, then for every $K\in \mathrm{Spec}^{%
\mathrm{fc}}(M)={\mathcal{S}}(M)$ we have $\{K\}={\mathcal{X}}(\{K\}_{e})$,
i.e. an open set. Since every singleton set is open, $\mathrm{Spec}^{\mathrm{%
fc}}(M)$ is discrete.

\item ($\Rightarrow $) Assume that $_{R}M$ has the a unique simple $R$%
-submodule. Clearly, $_{R}M$ has the min-property and $\mathrm{Spec}^{%
\mathrm{fc}}(M)$ is connected since it consists of only one point.

($\Leftarrow $) Assume that $_{R}M$ has the min-property and that $\mathrm{%
Spec}^{\mathrm{fc}}(M)$ is connected. By "1", $\mathrm{Spec}^{\mathrm{fc}%
}(M) $ is discrete and so ${\mathcal{S}}(M)=\mathrm{Spec}^{\mathrm{fc}}(M)$
has only one point since a discrete connected space cannot contain more than
one-point.$\blacksquare $
\end{enumerate}
\end{Beweis}

\begin{theorem}
\label{count-compact} Let $_{R}M$ be a S-IAD-module and assume that every
fully coprime $R$-submodule of $M$ is simple.

\begin{enumerate}
\item $\mathrm{Spec}^{\mathrm{fc}}(M)$ is countably compact if and only if ${%
\mathcal{S}}(M)$ is countable.

\item $\mathrm{Spec}^{\mathrm{fc}}(M)$ is compact if and only if ${\mathcal{S%
}}(M)$ is finite.
\end{enumerate}
\end{theorem}

\qquad As a direct consequence of Theorem \ref{compact} and Proposition \ref%
{it-irr} we obtain:

\begin{theorem}
\label{colocal} Let $_{R}M$ be a S-IAD-module and assume that every fully
coprime $R$-submodule of $M$ is simple. Then $_{R}M$ is colocal if and only
if $\mathrm{Spec}^{\mathrm{fc}}(M)$ is connected.
\end{theorem}


\begin{lemma}
\label{bireg}Let $_{R}M$ self-injective self-cogenerator duo and $S$ be
Noetherian with every prime ideal maximal. Then $\mathcal{S}(M)=\mathrm{Spec}%
^{\mathrm{fc}}(M)$.
\end{lemma}

\begin{Beweis}
Notice that $\mathcal{S}(M)\subseteq \mathrm{Spec}^{\mathrm{fc}}(M)$ (see
Remark \ref{sub-cop}). If $K\in \mathrm{Spec}^{\mathrm{fc}}(M),$ then $%
\mathrm{An}(K)\leq S$ is a prime ideal by Proposition \ref{corad=}, whence a
maximal ideal by our assumption on $S.$ It follows that $K=\mathrm{Ke}(%
\mathrm{An}(K))$ is simple: if $0\neq K_{1}\subsetneqq K,$ for some $%
K_{1}\leq _{R}M,$ then $\mathrm{An}(K)\subsetneqq \mathrm{An}%
(K_{1})\subsetneqq $ $S$ since $\mathrm{Ke}(-)$ is injective, a
contradiction.$\blacksquare $
\end{Beweis}

\begin{remark}
The ring $R$ is called $\pi$-\textit{regular}, iff for each $a \in R$ there
exist a positive integer $n = n(a)$, depending on $a$, and $x \in R$ such
that $a^n = a^n x a^n$. If $R$ is a left (right) duo ring, then every prime
ideal of $R$ is maximal if and only if $R$ is $\pi$-regular \cite{Hir1978}.
Moreover, $R$ is called \textit{biregular}, iff every ideal of $R$ is
generated by a central idempotent. By \cite[3.18 (6, 7)]{Wis1991}, every
prime ideal of a biregular ring is maximal.
\end{remark}

\begin{lemma}
\label{1n}Let $_{R}M$ be a top $R$-module. If $n\geq 2$ and $\mathcal{A}%
=\{K_{1},...,K_{n}\}\subseteq \mathrm{Spec}^{\mathrm{fc}}(M)$ is a connected
subset, then for every $i\in \{1,...,n\},$ there exists $j\in
\{1,...,n\}\backslash \{i\}$ such that $K_{i}\leq _{R}K_{j}$ or $K_{j}\leq
_{R}K_{i}.$
\end{lemma}

\begin{Beweis}
Without loss of generality, suppose $K_{1}\nsubseteqq K_{j}$ and $%
K_{j}\nsubseteqq K_{1}$ for all $2\leq j\leq n$ and set $F:=\dsum%
\limits_{i=2}^{n}K_{i},$ $W_{1}:=\mathcal{A}\cap \mathcal{X}^{\mathrm{fc}%
}(K_{1})=\{K_{2},...,K_{n}\}$ and $W_{2}:=\mathcal{A}\cap \mathcal{X}^{%
\mathrm{fc}}(F)=\{K_{1}\}$ (if $n=2,$ then clearly $W_{2}=\{K_{1}\};$ if $%
n>2 $ and $K_{1}\notin W_{2},$ then $K_{1}\subseteq
\dsum\limits_{i=2}^{n}K_{i}\subseteq
(K_{2}\odot_{M}\dsum\limits_{i=3}^{n}K_{i})$ and it follows that $%
K_{1}\subseteq \dsum\limits_{i=3}^{n}K_{i}.$ One can show by induction that $%
K_{1}\leq _{R}K_{n},$ a contradiction). So $\mathcal{A}=W_{1}\cup W_{2},$ a
disjoint union of proper non-empty open subsets, a contradiction.$%
\blacksquare $
\end{Beweis}

\begin{proposition}
\label{lf}Let $_{R}M$ be a S-IAD-module and let $\varnothing \neq \mathcal{K}%
=\{K_{\lambda }\}_{\Lambda }\subseteq \mathcal{S}(M).$ If $\left\vert
\mathcal{S}(L)\right\vert <\infty $ for every $L\in \mathrm{Spec}^{\mathrm{fc%
}}(M),$ then $\mathcal{K}$ is locally finite.
\end{proposition}

\begin{Beweis}
Let $L\in \mathrm{Spec}^{\mathrm{fc}}(M)$ and set
\begin{equation*}
F:=\sum\limits_{K\in \mathcal{K}\cap \mathcal{X}^{\mathrm{fc}}(L)}K\text{ \
\ (}:=0,\text{ if }\mathcal{K}\cap \mathcal{X}^{\mathrm{fc}}(L)=\varnothing
\text{)}
\end{equation*}%
Notice first that $L\nsubseteqq F:$ If $L\subseteq F,$ then there exists a
simple $R$-submodule $0\neq \widetilde{K}\subseteq L\subseteq F.$ Since $M$
has the min-property by Lemma \ref{min}, we conclude that $\widetilde{K}=K$
for some $K\in \mathcal{K}\cap \mathcal{X}^{\mathrm{fc}}(L)$, a
contradiction. So, $L\in \mathcal{X}^{\mathrm{fc}}(F).$ Since $\left\vert
\mathcal{S}(L)\right\vert <\infty ,$ there exists (if any) a finite number
of simple $R$-submodules $\{K_{\lambda _{1}},..,K_{\lambda _{n}}\}=\mathcal{K%
}\cap \mathcal{V}^{\mathrm{fc}}(L).$ It is clear that $\{K_{\lambda
_{1}},..,K_{\lambda _{n}}\}=\mathcal{K}\cap \mathcal{X}^{\mathrm{fc}}(F)$
and we are done.$\blacksquare $
\end{Beweis}

\begin{lemma}
If $_R M$ be a S-IAD-module, then the following are equivalent for any $L
\leq_R M$:

\begin{enumerate}
\item $L \in {\mathcal{S}}(M)$;

\item $L$ is fully coprime in $M$ and ${\mathcal{V}}^{\mathrm{fc}}(L) = \{
L\}$;

\item $\{L\}$ is closed in ${\mathbf{Z}}_M^{\mathrm{fc}}$.
\end{enumerate}
\end{lemma}

\begin{Beweis}
Notice that, by Remark \ref{sub-cop}, $\mathcal{S}(M)\subseteq \mathrm{Spec}%
^{\mathrm{fc}}(M)$. Let $L \leq_R M$.

$(1)\Rightarrow (2)$ and $(2)\Rightarrow (3)$ are obvious.

$(3)\Rightarrow (1)$ Assume that $\{L\}$ is closed in ${\mathbf{Z}}_{M}^{%
\mathrm{fc}}$. Then $\{L\}={\mathcal{V}}^{\mathrm{fc}}(K)$ for some $K\leq
_{R}M$. If $_{R}L$ is not simple then, since $_{R}M$ is atomic, there exists
some simple $R$-submodule $\tilde{L}$ of $M$ such that $\tilde{L}%
\varsubsetneqq L$. It follows that $\{L,\tilde{L}\}\subseteq {\mathcal{V}}^{%
\mathrm{fc}}(K)=\{L\}$, a contradiction.$\blacksquare $
\end{Beweis}

\qquad A topolgoical space is $T_{1}$ if and only if every singeleton set is
close. In light of this, the previous lemma yields:

\begin{proposition}
\label{T1} If $_{R}M$ is an S-IAD-module, then $\mathrm{Spec}^{\mathrm{fc}%
}(M)={\mathcal{S}}(M)$ if and only if $\mathbf{Z}^{\mathrm{fc}}(M)$ is $T_{1}
$ \emph{(}Fr\'{e}cht space\emph{)}.
\end{proposition}

\qquad Combining the prvious resutls we obtain

\begin{theorem}
\label{T2} Let $_{R}M$ be a S-IAD-module. The following are equivalent:

\begin{enumerate}
\item $\mathrm{Spec}^{\mathrm{fc}}(M)=\mathcal{S}(M);$

\item $\mathbf{Z}^{\mathrm{fc}}(M)$ is discrete;

\item $\mathbf{Z}^{\mathrm{fc}}(M)$ is $T_{2}$ \emph{(}Hausdorff space\emph{)%
};

\item $\mathbf{Z}^{\mathrm{fc}}(M)$ is $T_{1}$ \emph{(}Fr\'{e}cht space\emph{%
)}.
\end{enumerate}
\end{theorem}

\textbf{Acknowledgement}: The author thanks Professor Patrick Smith for
fruitful discussions on the topic during his visit to the University of
Glasgow (September 2008) and thereafter. He also thanks the University of
Glasgow for the hospitality and the Deanship of Scientific Research (DSR) at
King Fahd University of Petroleum $\&$ Minerals (KFUPM) for the financial
support.


\begin{thebibliography}{RRRF-AS2002}
\bibitem[Abu]{Abu} J. Y. Abuhlail, \emph{A Zariski topology for modules},
preprint $\langle $arXiv: math.RA 1007.3149$\rangle $.

\bibitem[Abu2008]{Abu2008} J. Y. Abuhlail, \emph{A Zariski topology for
bicomodules and corings, }Appl. Categ. Structures 16 (1-2) (2008), 13-28.

\bibitem[Abu2006]{Abu2006} J. Y. Abuhlail, \emph{Fully coprime comodules and
fully coprime corings, }Appl. Categ. Structures 14 (5-6) (2006), 379-409.

\bibitem[AS]{AS} Y. Al-Shaniafi and P. F. Smith, \emph{Comultiplication
modules over commutative rings}, to appear in Comm. Algebra.

\bibitem[AF1974]{AF1974} F. W. Anderson and K. Fuller, \emph{Rings and
Categories of Modules, }Springer-Verlag (1974).

\bibitem[AS2004]{AS2004} M. Ali and D. J. Smith, \emph{Some remarks on
multiplication and projective modules}, Comm. Algebra 32 (1) (2004),
3897--3909.

\bibitem[AM1969]{AM1969} M. Atiyah and I. Macdonald, \emph{Introduction to
commutative algebra,} Addison-Wesley Publishing Co. (1969).

\bibitem[A-TF2007]{A-TF2007} H. Ansari-Toroghy and H. Farshadifar, \emph{The
dual notion of multiplication modules, } Taiwanese J. Math. 11 (4) (2007),
1189--1201.

\bibitem[BJKN80]{BJKN80} L. Bican, P. Jambor, T. Kepka, P. N\u{e}mec, \emph{%
Prime and coprime modules, }Fund. Math. 107 (1) (1980)\textbf{,} 33-45.

\bibitem[Bou1998]{Bou1998} N. Bourbaki, \emph{Commutative algebra,}
Springer-Verlag (1998).

\bibitem[Bou1966]{Bou1966} N. Bourbaki, \emph{General Topology, Part I, }%
Addison-Wesley (1966).

\bibitem[CLVW2006]{CLVW2006} J. Clark, C. Lomp, N. Vanaja and R. Wisbauer,
\emph{Lifting Modules}, Birkh\"auser Verlag, Basel, 2006.

\bibitem[Dur1994]{Dur1994} T. Duraivel, \emph{Topology on spectrum of
modules, }J. Ramanujan Math. Soc. 9 (1) (1994), 25-34.

\bibitem[Dau1978]{Dau1978} J. Dauns, \emph{Prime modules, }J. Rein. Ang.
Math. 298 (1978), 165-181.

\bibitem[Hir1978]{Hir1978} Y. Hirano, \emph{Some studies on strongly $\pi$%
-regular rings}, Math. J. Okayama Univ. 20 (1978), 141--149.

\bibitem[HS2010]{HS2010} V. A. Hiremath and Poonam M. Shanbhag, \emph{Atomic
modules, }Int. J. Algebra, 4 (2) (2010), 61-69.

\bibitem[Joh53]{Joh53} R. Johnson, \emph{Representations of prime rings},
Trans. Amer. Math. Soc. 74 (1953), 351-357.

\bibitem[Lu1984]{Lu1984} Chin-Pi Lu, \emph{Prime submodules of modules},
Comment. Math. Univ. St. Paul. 33(1) (1984)\textbf{,} 61-69.

\bibitem[Lu1999]{Lu1999} C-P. Lu, \emph{The Zariski topology on the prime
spectrum of a module, }Houston J. Math. 25 (3) (1999), 417-432.

\bibitem[LY2006]{LY2006} D. Lu and W. Yu, \emph{On prime spectrum of
commutative rings}, Comm. Algebra 34 (2006), 2667--2672.

\bibitem[MMS1997]{MMS1997} R. McCasland, M. Moore and P. Smith, \emph{On the
spectrum of a module over a commutative ring}. Comm. Algebra 25 (1997),
79-103.

\bibitem[MMS1998]{MMS1998} R. McCasland, M. Moore and P. Smith, \emph{An
introduction to Zariski spaces over Zariski topologies}, Rocky Mountain J.
Math. 28(4) (1998), 1357-1369.

\bibitem[NT2001]{NT2001} R. Nekooei and L. Torkzadeh, \emph{Topology on
coalgebras, }Bull. Iran. Math. Soc. 27(2) (2001), 45-63.

\bibitem[OHS2006]{OHS2006} A. \"{O}zcan, A. Harmanci and P. F. Smith, \emph{%
Duo modules}, Glasgow Math. J. 48 (2006), 33---545.

\bibitem[PC1995]{PC1995} Y. S. Park and C. W. Choi, \emph{Multiplication
modules and characteristic submodules}, Bull. Korean Math. Soc. 32 (2)
(1995), 321--327.

\bibitem[RRRF-AS2002]{RRRF-AS2002} F. Raggi, J. Rios, H. Rinc\'{o}n, R. Fern%
\'{a}ndez-Alonso, C. Signoret, \emph{The lattice structure of preradicals.
II, Partitions}, J. Algebra Appl. 1(2) (2002), 201-214.

\bibitem[RRRF-AS2005]{RRRF-AS2005} F. Raggi, J. R\'{\i}os, H. Rinc\'{o}n, R/
Fern\`{a}ndez-Alonso and C. Signoret, \emph{Prime and irreducible preradicals%
}, J. Alg. Appl. 4 (4) (2005), 451--466.

\bibitem[RRW2005]{RRW2005} F. Raggi, J. R\'{\i}os Montes and R. Wisbauer,
\emph{Coprime preradicals and modules, }J. Pur. App. Alg. 200 (2005), 51-69.

\bibitem[Smi]{Smi} P. Smith, \emph{Modules with coindependent maximal
submodules}, to appear in J. Alg. Appl.

\bibitem[Smi1994]{Smi1994} P. Smith, \emph{Multiplication modules and
projective modules}, Period. Math. Hungar. 29 (2) (1994), 163--168.

\bibitem[ST2010]{ST2010} N. Schwartz and M. Tressl, \emph{Elementary
properties of minimal and maximal points in Zariski spectra}, J. Algebra 323
(2010), 698--728.

\bibitem[Tug2003]{Tug2003} A. A. Tuganbaev, \emph{\ Multiplication modules
over noncommutative rings}, Sb. Math. 194 (11-12) (2003), 1837--1864.

\bibitem[Tug2004]{Tug2004} A. A. Tuganbaev, \emph{Multiplication modules},
J. Math. Sci. (N.Y.) 123(2) (2004), 3839-3905.

\bibitem[Wij2006]{Wij2006} I. Wijayanti, \emph{Coprime modules and
comodules, }Ph.D. Dissertation, Heinrich-Heine Universit\"{a}t, D\"{u}%
sseldorf (2006).

\bibitem[Wis1991]{Wis1991} R. Wisbauer, \emph{Foundations of module and ring
theory. A handbook for study and research.} Gordon and Breach Science
Publishers (1991).

\bibitem[Wis1996]{Wis1996} R. Wisbauer, \emph{Modules and algebras :
Bimodule structure and\ group action on algebras}, Addison Wesely Longman
Limited (1996).

\bibitem[WW2009]{WW2009} I. Wijayanti and R. Wisbauer, \emph{Coprime modules
and comodules, }Commun. Algebra 37 (4) (2009), 1308--1333.

\bibitem[Zha2006-a]{Zha2006-a} Guoyin Zhang, \emph{Properties of top modules}%
, Int. J. Pure Appl. Math. 31 (3) (2006), 297--306.

\bibitem[Zha2006-b]{Zha2006-b} Guoyin Zhang, \emph{Multiplication modules
over any rings}, Nanjing Daxue Xuebao Shuxue Bannian Kan 23 (1) (2006),
59--69.

\bibitem[ZTW2006]{ZTW2006} G. Zhang, W. Tong and F. Wang, \emph{Spectrum of
a noncommutative ring}, Comm. Algebra 34 (8) (2006), 2795--2810.
\end{thebibliography}
\end{document}